\documentclass[
 sd,
amsmath,amssymb,
preprint,%
$reprint,%
numerical,
]{revtex4-1}
\usepackage{graphicx}
\usepackage{dcolumn}
\usepackage{bm}
\usepackage{hyperref}
\usepackage[utf8]{inputenc}
\usepackage[T1]{fontenc}
\usepackage{mathptmx}
\usepackage{etoolbox}
\makeatletter
\def\@email#1#2{%
	\endgroup
	\patchcmd{\titleblock@produce}
	{\frontmatter@RRAPformat}
	{\frontmatter@RRAPformat{\produce@RRAP{*#1\href{mailto:#2}{#2}}}\frontmatter@RRAPformat}
	{}{}
}%
\makeatother
\begin{document}

	\title{A combined effect of rigid top surface with diffuse and oblique collimated irradiation on the stability of the suspension of phototactic microorganisms }
	\author{S. K. Rajput}
	\email{shubh.iiitj@gmail.com}
	\affiliation{ 
		Department of Mathematics, PDPM Indian Institute of Information Technology Design and Manufacturing,
		Jabalpur 482005, India.
	}%
	
	

	\begin{abstract}
		This article explores how a rigid top surface with diffuse and oblique collimated irradiation affect isotropic scattering algal suspensions. When the fluid flow becomes zero, the suspension reaches a steady (basic) state where up-and-down swimming occurs due to interplay of phototaxis and diffusion. For purely scattering suspensions, a bimodal steady state occurs due to scattering, which reverts to a unimodal steady state as the angle of incidence increases with fixed other governing parameters. To check the linear stability of the suspension, a small perturbation to the basic state are considered and the perturbed equations are solved by using the Newton-Raphson-Kantorovich (NRK) iterative method. The linear stability of the same suspension predicts both stable and oscillatory nature of disturbance for specific parameter ranges. A rigid top surface, as well as diffuse and oblique collimated irradiation, make the suspension more stable.
		
	\end{abstract}
	
	\maketitle

	\section{INTRODUCTION}
	
	Bio-convection refers to the collective motion of microorganisms or biological components, such as motile algae, bacteria, or spermatozoa, that can generate fluid flow patterns in the medium (water here). These microorganisms are denser than the medium and they swim upward on average. Pattern formation depends on the swimming behavior of the microorganisms. Whenever microorganisms stop swimming, the pattern formation in bio-convection disappears. However, up-swimming and higher density are not necessary for pattern formation. On the other hand, the pattern formation (or swimming behavior of microorganisms) is affected by taxis which refers to the directed movement of motile organisms or cells in response to a environmental stimulus, such as light, chemicals, or gravity. The stimulus can be attractive or repulsive, and the movement can be towards or away from the stimulus. Taxis is an essential mechanism for many organisms to find food, avoid predators, and navigate their environment. For example, phototaxis is the movement of an organism towards the light (positive phototaxis) or away from light (negative phototaxis), while chemotaxis is the movement towards or away from a chemical gradient. There are other examples of taxis such as gravitaxis and gyrotaxis etc. Here, this article related to the phototaxis only.
	 
	 The pattern formation in bioconvection may be significantly influenced by various forms of illumination intensity (such as diffuse irradiation)~\cite{1wager1911,2kitsunezaki2007}. 
	Strong (bright) light can damage the stable patterns or prevents the formation of patterns in a suspension of motile microorganisms.
	The pattern's size, shape, and structure may all be affected by the light intensity~\cite{3kessler1985,4williams2011,5kessler1989}. Variations in bioconvection patterns caused by light intensity can be explained through the following circumstances. First, the phototactic nature of the microorganisms. The cells move towards the light source for obtain energy via photosynthesis but they change their direction to avoid photo damage. Thus, the algae cells try to accumulate at a suitable location where they can find optimal light intensity. Second, the self-shading and scattering of light may be affect the pattern formation~\cite{7ghorai2010}. Finally, in the presence of the diffuse irradiation, the patterns can be affected due to the uniformity of the diffuse irradiation.\par
	
	 we utilize the phototaxis model developed by Panda $et$ $al$.~\cite{15panda2016}. the diffuse solar radiation is one of the part of sunlight which occurs due to scattering of the direct light in the presence of the clouds. Therefore, to compensate the solar loses, the algal suspension can be consider illuminated by both diffuse and oblique irradiation and advantage of the diffuse irradiation can not be ignored in building effective photo-bioreactors. The phototactic algae is used for fixation of carbon dioxide in photo-bioreactors, and biofuels may also made from the produced biomass. In the biofuel production, phototactic bio-convection can be less worthy. Therefore,  to understand of behavior of algal species in the suspension, a realistic phototaxis model should consider the effects of bothe diffuse and collimated  irradiation.\par
	
	Consider the bioconvection in a dilute suspension of phototactic algae. The basic state (sub layer) of same suspension is formed where cells motion are driven by the interplay of phototaxis and diffusion of cells. This sublayer can form at different depths in the suspension, depending on the total intensity of light and the critical intensity $G_c$, which is the value of the total intensity at which the sublayer forms at a particular depth. When the total intensity is lower than $G_c$ everywhere in the suspension, the sublayer forms at the top of the suspension. Conversely, if the total intensity is higher than $G_c$ everywhere, the sublayer forms at the bottom. When the total intensity is equal to $G_c$ at an interior depth of the suspension, the sublayer forms approximately at that depth. The region below the sublayer is gravitationally unstable, while the region above is stable. This means that if the bioconvection system becomes unstable, the fluid motions in the lower unstable region will penetrate into the upper stable region. This phenomenon is known as penetrative convection and is observed in many other convection problems as well~\cite{9straughan1993,10ghorai2005,11panda2016}.
	
	\begin{figure}[!h]
		\centering
		\includegraphics[width=14cm]{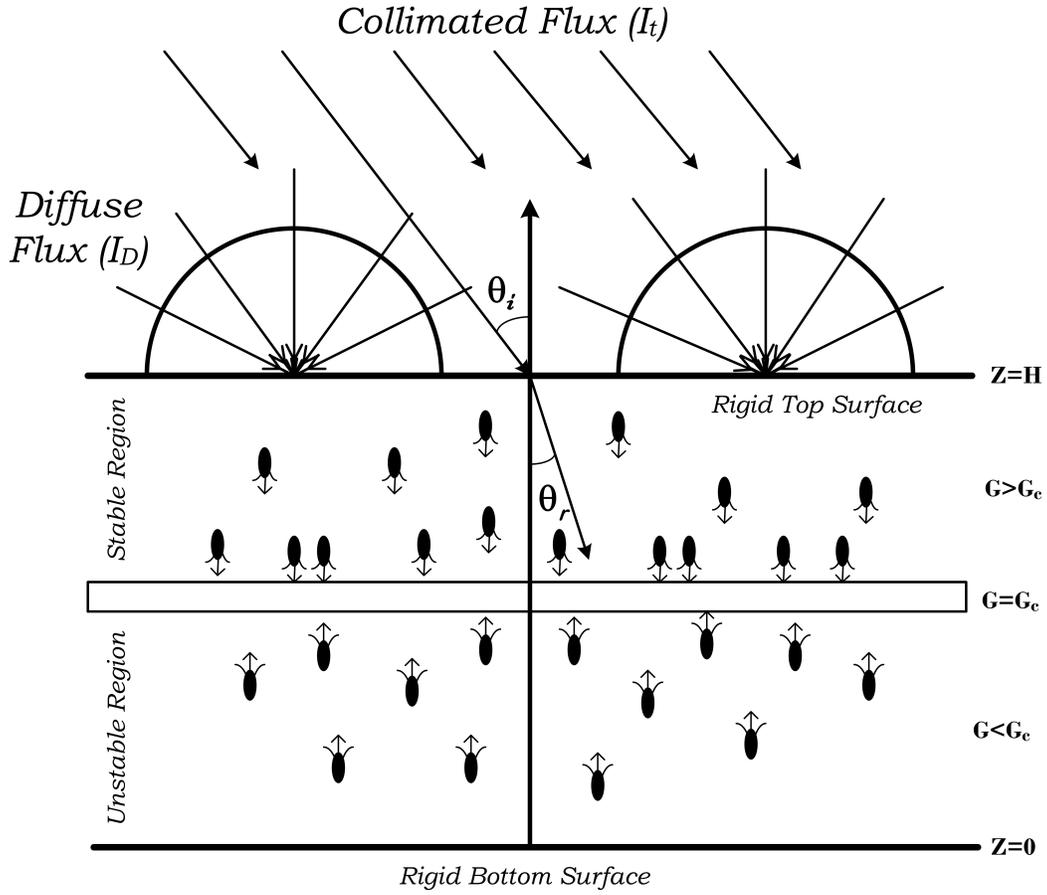}
		\caption{\footnotesize{Formation of the sublayer in the interior of the algal suspension at $G=G_c$, where $G_c$ is the critical light intensity. Above (below) the sublayer, the suspension is stable (unstable).}}
		\label{fig1}
	\end{figure}
	
	The various studies have investigated the phenomenon of phototactic bioconvection in suspensions of phototactic algae under different conditions. Vincent and Hill~\cite{12vincent1996} explored bioconvection in a phototactic algae suspension, analyzing the equilibrium solution and finding stationary and oscillatory modes of disturbance at the onset of bioconvective instability. Ghorai and Hill~\cite{10ghorai2005} studied two-dimensional phototactic bioconvection numerically using the Vincent and Hill model, but did not account for scattering effects. Ghorai et al.~\cite{7ghorai2010} examined the onset of bioconvection in a suspension of isotropically scattering phototactic algae and found an unusual bimodal steady-state profile and oscillatory instabilities for certain parameter values. Ghorai and Panda~\cite{14panda2013} cheached  the stability of anisotropic scattering algal suspension, observing a transition from stationary to oscillatory modes with the variation of anisotropic scattering coefficient for certain parameter values. Panda and Ghorai~\cite{14panda2013} simulated two-dimensional phototactic bioconvection in a non-linear regime of an absorbing and isotropic scattering suspension, with results differing from those found by Ghorai and Hill~\cite{10ghorai2005} due to scattering effects. Panda and Singh~\cite{11panda2016} investigated two-dimensional phototactic bioconvection in a suspension confined by lateral/side walls, observing a significant stabilizing effect on the suspension due to the presence of walls. Panda et al.~\cite{15panda2016} studied the effects of diffuse irradiation on an isotropic light scattering suspension of phototactic algae and observed a stabilizing effect due to diffuse irradiation, with a transition from bimodal to unimodal vertical concentration profiles at the base state.Panda~\cite{8panda2020} studied the impact of forward anisotropic scattering on the onset of phototactic bioconvection using both diffuse and collimated irradiation. Panda $et$ $al$.~\cite{16panda2022} studied the effect of oblique irradiation on algal suspensions and observed that the location of the maximum concentration of microorganisms in a suspension shifts towards the top of the suspension and the value of the maximum concentration increases as the angle of incidence increases. In the same article, they checked the linear stability of the suspension for different angles of incidence. Recently, Kumar~\cite{17kumar2022} investigated the effect of oblique collimated irradiation on the isotropic scattering algal suspension. He found both types of solutions (stationary and overstable) for certain ranges of parameters. More recently, Kumar~\cite{39kumar2023} studied the effect of collimated irradiation on the algae suspension where both vertical walls were assumed to be rigid. In their study, he found the stabilizing effect on the suspension due to rigid walls. However, no study to date has explored the onset of phototactic bioconvection that incorporates the effects of both diffuse and oblique collimated irradiation with rigid upper surface. Therefore, this study investigates the effects of rigid top surface of the suspension which is illuminated by both types of irradiation.\par
	
The article is structured as follows: Firstly, the problem is mathematically formulated. The equilibrium solution is then obtained and the base bioconvective governing system is perturbed by small disturbances. The linear stability problem is derived and solved numerically. The results of the model are then presented and compared to a non-scattering phototaxis model. Finally, the novelty of the proposed model is discussed.
	
	\section{MATHEMATICAL FORMULATION}
	
		In this model, the motion is considered in a dilute suspension of phototactic algae in a layer of limited depth $H$, but infinite width, under illumination from oblique collimated and diffuse irradiation from above. The top and bottom boundaries of the layer are assumed to be non-reflective. The intensity of light at a particular location $x$ in a unit direction $\boldsymbol{ s}$ is represented by $I(x,\boldsymbol{ s})$, where $x$ is measured relative to a rectangular Cartesian coordinate system $O(x,y,z)$ with the $z$-axis oriented vertically upward, and $s$ is defined by the angle $\theta$ with respect to the $z$-axis and the angles $\phi$ and $\theta$ with respect to the $x$ and $y$-axes, respectively.
	 

	\subsection{\label{sec:level3}PHOTOTAXIS WITH ABSORPTION AND SCATTERING}
	
	The RTE is utilized to govern the light intensity in the absorbing and scattering medium, which is given by 
	\begin{equation}\label{1}
		\frac{dI(\boldsymbol{x},\boldsymbol{s})}{ds}+(\alpha+\beta_s)I(\boldsymbol{x},\boldsymbol{s})=\frac{\beta_s}{4\pi}\int_{0}^{4\pi}I(\boldsymbol{x},\boldsymbol{s'})\varphi(\boldsymbol{s},\boldsymbol{s'})d\Omega',
	\end{equation}
	where $\alpha$ and $\beta_s$ are the absorption and scattering coefficients, respectively, and $\varphi(\boldsymbol{s},\boldsymbol{s'})$ is the scattering phase function. In this model, isotropic scattering is consider for simplicity. Therefore, we use $\varphi(\boldsymbol{s},\boldsymbol{s'})=1$ here ~\cite{15panda2016}.

	\begin{figure}[!h]
		\centering
		\includegraphics[width=14cm]{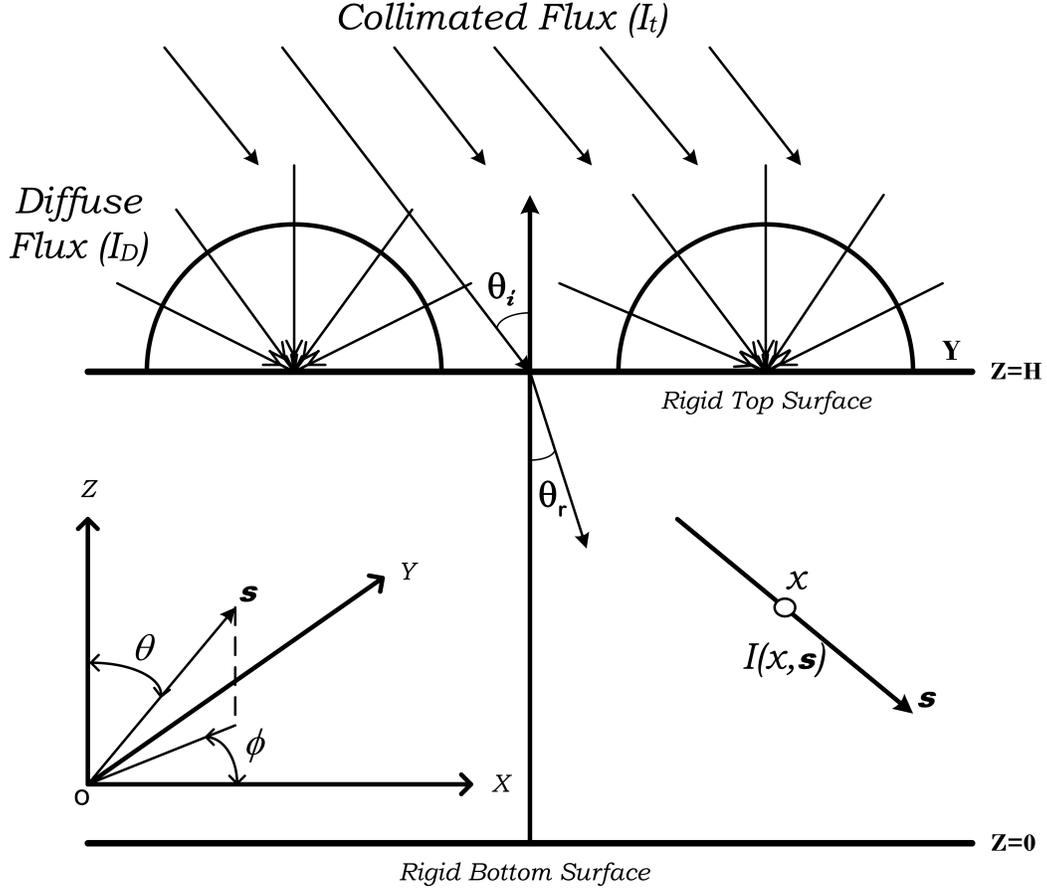}
		\caption{\footnotesize{Geometric configuration of the problem.}}
		\label{fig2}
	\end{figure}

	The light intensity on the top of the suspension at $\boldsymbol{x}_H=(x,y,H)$ is given by
	\begin{equation*}
		I(\boldsymbol{x}_H,\boldsymbol{s})=I_t\delta(\boldsymbol{s}-\boldsymbol{s_0})+\frac{I_D}{\pi}, 
	\end{equation*}
	where $I_t$ and $I_D$ are the magnitudes of direct (collimated) and diffuse irradiation respectively~\cite{8panda2020,15panda2016}. Consider, the light absorption and scattering is proportional to the number of cells. Therefore, $\alpha=a n(\boldsymbol{x})$ and $\beta_s=b n(\boldsymbol{ x})$, then the RTE becomes
	\begin{equation}\label{2}
		\frac{dI(\boldsymbol{x},\boldsymbol{s})}{ds}+(a+a)nI(\boldsymbol{x},\boldsymbol{s})=\frac{b n}{4\pi}\int_{0}^{4\pi}I(\boldsymbol{x},\boldsymbol{s'})d\Omega'.
	\end{equation}
	The total intensity at a fixed point $\boldsymbol{x}$ in the medium is 
	\begin{equation*}
		G(\boldsymbol{x})=\int_0^{4\pi}I(\boldsymbol{x},\boldsymbol{s})d\Omega,
	\end{equation*}
	and the radiative heat flux is defined as
	\begin{equation}\label{3}
		\boldsymbol{q}(\boldsymbol{x})=\int_0^{4\pi}I(\boldsymbol{x},\boldsymbol{s})\boldsymbol{s}d\Omega.
	\end{equation}
	Let $<\boldsymbol{\wp}>$ be the average swimming direction where $\boldsymbol{ \wp}$ is the unit vector in the direction of cell swimming. The swimming speed of microorganism  is independent from the illumination, position, time, and direction for various species of microorganisms~\cite{18hill1997}. In this model, the swimming speed of cells and the fluid is consider same, then the mean cell
	swimming velocity is given by
	\begin{equation*}
		\boldsymbol{U}_c=U_c<\boldsymbol{p}>,
	\end{equation*}
	where $U_c$ is the average cell swimming speed and the mean cell swimming direction $<\boldsymbol{\wp}>$ is given by
	\begin{equation}\label{4}
		<\boldsymbol{\wp}>=-T(G)\frac{\boldsymbol{q}}{\varkappa+|\boldsymbol{q}|},
	\end{equation}
	where $T(G), $ is taxis function, which has the mathematical form as 
	\begin{equation*}
		T(G)=\left\{\begin{array}{ll}\geq 0, & \mbox{if } G(\boldsymbol{x})\leq G_{c}, \\
			< 0, & \mbox{if }G(\boldsymbol{x})>G_{c}.  \end{array}\right. 
	\end{equation*}
	
	The exact functional form of taxis function depends on the species of the microorganisms~\cite{12vincent1996}. The non-negative $\varkappa$ is introduced to handle the case of isotropic light conditions, but here the light intensity throughout the suspension is not isotropic. Therefore, $\varkappa=0$ is use here.~\cite{7ghorai2010}

	\subsection{GOVERNING EQUATIONS}
	
	In this model, we assume a monodisperse cell population in dilute suspension $(n\mathtt{v}<<1)$, therefore cell's volume is small and interactions among cells are negligible. The volume of each cell $\mathtt{v}$ and density $\rho+\Delta\rho$, where $\rho$ is the water density  and $\Delta\rho<<\rho$. The average fluid velocity is $\boldsymbol{u}$, and the concentration of algal cells is $n$ in the unite volume. In this model, we deal with incompressible suspension, then the equation of continuity is
\begin{equation}\label{5}
	\boldsymbol{\nabla}\cdot \boldsymbol{u}=0.
\end{equation}
We assume that stokeslets predominate the effect of concentration of cells on the suspension except for their negative buoyancy for simplicity.
Therefore, under the Boussinesq approximation the momentum equation is 
\begin{equation}\label{6}
	\rho\frac{D\boldsymbol{u}}{Dt}=-\boldsymbol{\nabla} P_e+\mu\nabla^2\boldsymbol{u}-nvg\Delta\rho\hat{\boldsymbol{z}},
\end{equation}
where $D/Dt=\partial/\partial t+\boldsymbol{u}\cdot\boldsymbol{\nabla}$ is the material, $P_e$ is the excess pressure and $\mu$ is the suspension's viscosity which is assumed to be that of fluid, here.\newline 
The cell conservation equation is given by
\begin{equation}\label{7}
	\frac{\partial n}{\partial t}=-\boldsymbol{\nabla}\cdot \boldsymbol{J},
\end{equation}
where $\boldsymbol{J}$ is the total cell flux which is given by
\begin{equation}\label{8}
	\boldsymbol{J}=n(u+W_c<\boldsymbol{\wp}>)-\boldsymbol{D}\boldsymbol{\nabla} n.
\end{equation}
On R.H.S. of Eq.~(\ref{8}) the first term is arises due to the advection of cells by the bulk fluid flow, the second term occurs due to the average swimming of the cells and the third term arises due to diffusion of cells. Here, the diffusivity tensor $\boldsymbol{D}$ is consider isotropic and constant such that $\boldsymbol{D} = DI$. The representation of cell flux in Eq.~(\ref{8}) has two main assumptions. First, cells are purely phototactic, so the effect of viscous torque is ignored and second, diffusion tensor has constant value. These two assumptions have a great importance in this model. With these assumptions, we can remove the Fokker-Planck equation from the governing equations and the resulting model can be used as a limiting case to estimate the difficulty of the problem before constructing a more complex model.

\subsection{BOUNDARY CONDITIONS}
Consider the lower and upper boundaries are rigid in this model.
Thus, the boundary conditions are
\begin{equation}\label{9}
	\boldsymbol{u}\times\hat{\boldsymbol{z}}=0\qquad on\quad z=0,H,
\end{equation}
\begin{equation}\label{10}
	\boldsymbol{J}\cdot\hat{\boldsymbol{z}}=0\qquad on\quad z=0,H.
\end{equation}

We assume that the suspension is uniformly illuminated by both diffuse  and collimated oblique irradiation, then the boundary condition for intensities are
\begin{subequations}
	\begin{equation}\label{11a}
		I(x, y, z = 1, \theta, \phi)=I_t\delta(\boldsymbol{s}-\boldsymbol{s_0})+\frac{I_D}{\pi},\quad (\pi/2\leq\theta\leq\pi),
	\end{equation}
	\begin{equation}\label{11b}
		I(x, y, z = 0, \theta, \phi) =0,\quad (0\leq\theta\leq\pi/2).
	\end{equation}
\end{subequations}

\subsection{DIMENSIONLESS EQUATIONS}
The governing equations are made dimensionless by scaling all lengths on $H$, the depth of the
layer, time on the diffusive time scale $H^2/D$, and the bulk fluid velocity on $D/H$. The appropriate
scaling for the pressure is $\mu D/H^2$, and the cell concentration is scaled on $\bar{n}$, the mean concentration.
In terms of the non-dimensional variables, the governing equations become
\begin{equation}\label{12}
	\boldsymbol{\nabla}\cdot\boldsymbol{u}=0,
\end{equation}
\begin{equation}\label{13}
	S_{c}^{-1}\left(\frac{D\boldsymbol{u}}{Dt}\right)=-\nabla P_{e}-Rn\hat{\boldsymbol{z}}+\nabla^{2}\boldsymbol{u},
\end{equation}
\begin{equation}\label{14}
	\frac{\partial{n}}{\partial{t}}=-{\boldsymbol{\nabla}}\cdot{\boldsymbol{J}},
\end{equation}
where
\begin{equation}\label{15}
	{\boldsymbol{J}}=n{\boldsymbol{u}}+nV_{c}<{\boldsymbol{\wp}}>-{\boldsymbol{\nabla}}n,
\end{equation}

where $S_{c}=\nu/{D}$ is the Schmidt number, $V_c=W_cH/D$ is scaled swimming speed, and $R=\bar{n}v g\Delta{\rho}H^{3}/\nu\rho{D}$ is the Rayleigh number.
In dimensionless form, the boundary conditions become
\begin{equation}\label{16}
	\boldsymbol{u}\times\hat{\boldsymbol{z}}=0\qquad on\quad z=0,1,
\end{equation}
\begin{equation}\label{17}
	\boldsymbol{J}\cdot\hat{\boldsymbol{z}}=0\qquad on\quad z=0,1.
\end{equation}

Then, in dimensionless form  RTE becomes
\begin{equation}\label{18}
	\frac{dI(\boldsymbol{x},\boldsymbol{s})}{ds}+\kappa nI(\boldsymbol{x},\boldsymbol{s})=\frac{\sigma_s n}{4\pi}\int_{0}^{4\pi}I(\boldsymbol{x},\boldsymbol{s'})d\Omega',
\end{equation}
where $\kappa=(a+b)\Bar{n}H$,and $\sigma_s=b\Bar{n}H$ are the extinction and scattering coefficient in dimensionless form ,respectively. The scattering albedo $\omega=\sigma_s/\kappa$ is used to measure the scattering efficiency of microorganisms, here. In terms of scattering albedo $\omega$, Eq.~(\ref{18}) becomes
\begin{equation}\label{19}
	\frac{dI(\boldsymbol{x},\boldsymbol{s})}{ds}+\kappa nI(\boldsymbol{x},\boldsymbol{s})=\frac{\omega\kappa n}{4\pi}\int_{0}^{4\pi}I(\boldsymbol{x},\boldsymbol{s'})d\Omega',
\end{equation}
where scattering albedo $\omega\in$ [0, 1]. The value of the scattering albedo $\omega=0$ is used for purely absorbing medium while $\omega=1$ implies purely scattering medium. RTE in the form of direction cosine is 

\begin{equation}\label{20}
	\xi\frac{dI}{dx}+\eta\frac{dI}{dy}+\nu\frac{dI}{dz}+\kappa nI(\boldsymbol{x},\boldsymbol{s})=\frac{\omega\kappa n}{4\pi}\int_{0}^{4\pi}I(\boldsymbol{x},\boldsymbol{s'})d\Omega',
\end{equation}
where $\xi=sin\theta\cos\phi,\eta=\sin\theta\sin\phi$ and $\nu=\cos\theta$ are the direction cosines in x, y and z direction. In dimensionless form, the intensity at boundaries becomes,
\begin{subequations}
	\begin{equation}\label{21a}
		I(x, y, z = 1, \theta, \phi)=I_t\delta(\boldsymbol{s}-\boldsymbol{s_0})+\frac{I_D}{\pi} ,\qquad (\pi/2\leq\theta\leq\pi),
	\end{equation}
	\begin{equation}\label{21b}
		I(x, y, z = 0, \theta, \phi) =0,\qquad (0\leq\theta\leq\pi/2). 
	\end{equation}
\end{subequations}
	
	\section{THE BASIC (EQUILIBRIUM) STATE SOLUTION}
	
	Equations $(\ref{12})-(\ref{15})$ and $(\ref{20})$ with the appropriate boundary conditions possess an equilibrium solution in which
	
	\begin{equation}\label{22}
		\boldsymbol{u}=0,~~~n=n_s(z)\quad and\quad  I=I_s(z,\theta).
	\end{equation}
	Therefore, in the basic state total intensity $G_s$ and radiative heat flux $\boldsymbol{q}_s$ are given by the relations
	
	\begin{equation*}
		G_s=\int_0^{4\pi}I_s(z,\theta)d\Omega,\quad 
		\boldsymbol{q}_s=\int_0^{4\pi}I_s(z,\theta)\boldsymbol{s}d\Omega.
	\end{equation*}
	 Intensity in the basic state, $I_s$, can be govern by the equation
	\begin{equation}\label{23}
		\frac{dI_s}{dz}+\frac{\kappa n_sI_s}{\nu}=\frac{\omega\kappa n_s}{4\pi\nu}G_s(z).
	\end{equation}
	
	The basic state intensity can be decomposed into collimated part $I_s^c$ and diffuse part $I_s^d$ such that $I_s=I_s^c+I_s^d$. The collimated part part of the basic state intensity $I_s^c$ is calculated by the equation 
	
	\begin{equation}\label{24}
		\frac{dI_s^c}{dz}+\frac{\kappa n_sI_s^c}{\nu}=0,
	\end{equation}
	
	with the boundary condition

	\begin{equation}\label{25}
		I_s^c( 1, \theta) =I_t\delta(\boldsymbol{s}-\boldsymbol{ s}_0),\qquad (\pi/2\leq\theta\leq\pi). 
	\end{equation}
	Therefore, we find $I_s^c$ in the form of 
	\begin{equation}\label{26}
		I_s^c=I_t\exp\left(\int_z^1\frac{\kappa n_s(z')}{\nu}dz'\right)\delta(\boldsymbol{s}-\boldsymbol{s_0}), 
	\end{equation}
	
	and the diffused part is calculated by  
	\begin{equation}\label{27}
		\frac{dI_s^d}{dz}+\frac{\kappa n_sI_s^d}{\nu}=\frac{\omega\kappa n_s}{4\pi\nu}G_s(z),
	\end{equation}
	with the boundary conditions
	\begin{subequations}
		\begin{equation}\label{28a}
			I_s^d( 1, \theta) =\frac{I_D}{\pi},\qquad (\pi/2\leq\theta\leq\pi), 
		\end{equation}
		\begin{equation}\label{28b}
			I_s^d( 0, \theta) =0,\qquad (0\leq\theta\leq\pi/2). 
		\end{equation}
	\end{subequations}
	
	Now the total intensity, $G_s=G_s^c+G_s^d$ in the equilibrium state can be written as
	\begin{equation}\label{29}
		G_s^c=\int_0^{4\pi}I_s^c(z,\theta)d\Omega=I_t\exp\left(\frac{-\int_z^1\kappa n_s(z')dz'}{\cos\theta_r}\right),
	\end{equation}
	\begin{equation}\label{30}
		G_s^d=\int_0^{\pi}I_s^d(z,\theta)d\Omega.
	\end{equation}
	
	If we define a new variable as 
	\begin{equation*}
		\tau=\int_z^1 \kappa n_s(z')dz',
	\end{equation*}
	
	then dimensionless total intensity, $\Lambda(\tau)=G_s(\tau)/I_t$, satisfies the following Fredholm Integral Equation (FIE),
	
	\begin{equation}\label{31}
		\Lambda(\tau) = \frac{\omega}{2}\int_0^\kappa \Lambda(\tau')E_1(|\tau-\tau'|)d\tau'+e^{-\tau/\cos\theta_r}+2I_DE_2(\tau),
	\end{equation}
	
	he,re $E_1(x)$ and $E_3(x)$ are the exponential integral of order 1 and 2, respectively. This FIE is solved by using the method of subtraction of singularity.\par
	
	The radiative heat flux in the basic state is written as
	
	\begin{equation*}
		\boldsymbol{q_s}=\int_0^{4\pi}\left(I_s^c(z,\theta)+I_s^d(z,\theta)\right)\boldsymbol{s}d\Omega=-I_t(\cos\theta_r)\exp\left(\frac{\int_z^1-\kappa n_s(z')dz'}{cos\theta_r}\right)\hat{\boldsymbol{z}}+\int_0^{4\pi}I_s^d(z,\theta)\boldsymbol{s}d\Omega.
	\end{equation*}
	
	since $I_s^d(z,\theta)$ is not dependent on $\phi$, so the horizontal components of $\boldsymbol{q_s}$ vanish. Therefore, in the basic state the radiative heat flux $\boldsymbol{q}_s=-q_s\hat{\boldsymbol{z}}$, where $q_s=|\boldsymbol{q_s}|$. Then the mean swimming direction is calculated by
	
	\begin{equation*}
		<\boldsymbol{\wp_s}>=-T_s\frac{\boldsymbol{q_s}}{q_s}=T_s\hat{\boldsymbol{z}},
	\end{equation*}
	
	where $T_s=T(G_s).$\par
	The basic concentration of algae cells $n_s(z)$ satisfies the following equation
	
	\begin{equation}\label{32}
		\frac{dn_s}{dz}-V_cT_sn_s=0,
	\end{equation}
	which is augmented by the cell conservation relation
	\begin{equation}\label{33}
		\int_0^1n_s(z)dz=1.
	\end{equation}
	Eqs.~(\ref{31})-(\ref{33}) represent a boundary value problem (BVP) and this BVP is solved by using a
	shooting method numerically.
	
	\begin{figure*}[!bt]
		\includegraphics{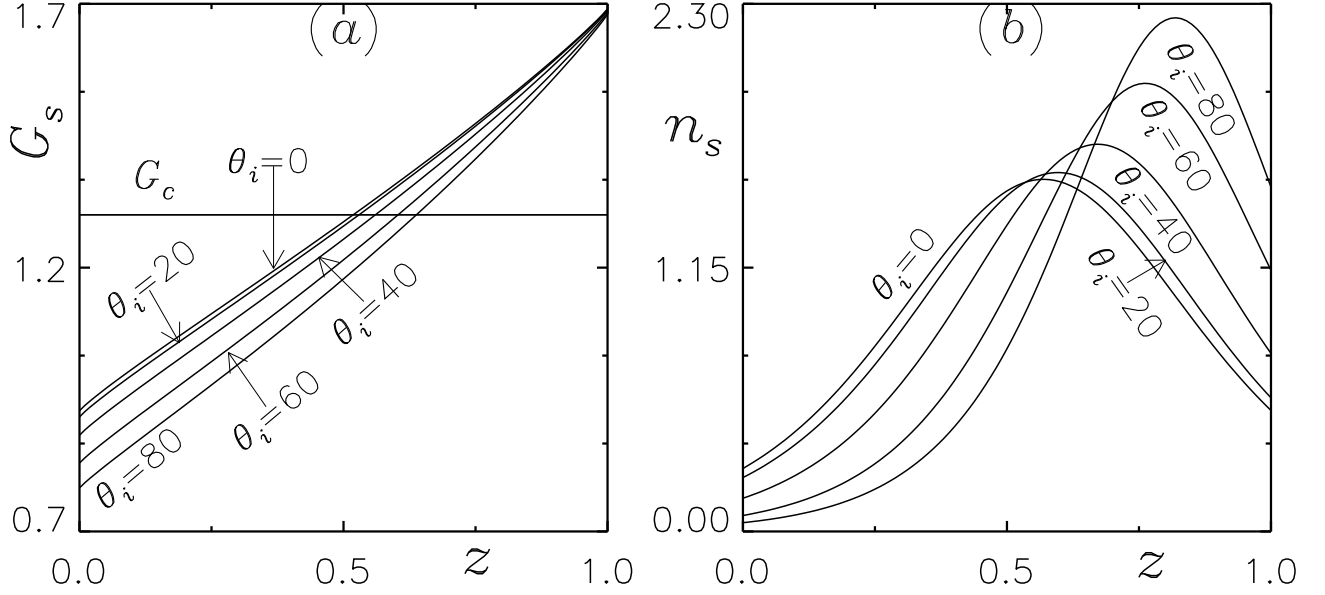}
		\caption{\label{fig3}(a) Variation of total intensity in a uniform suspension for discrete values of the angle of incidence $\theta_i$, (b) corresponding base concentration profile for the governing parameter values $S_c=20, V_c=10,k=0.5, I_D=0.25,\omega=0.4$ and $I_t=1$. Here, the critical intensity $G_c=1.3$ is utilized.}
	\end{figure*}

	The incident radiation intensity, $I_t=1$ on the top is considered in this article.   
	Consider the phototaxis function
	\begin{equation}\label{34}
		T(G)=0.8\sin\left(\frac{3\pi}{2}\chi(G)\right)-0.1\sin\left(\frac{\pi}{2}\chi(G)\right),\quad \chi(G)=\frac{G}{3.8}\exp[0.252(3.8-G)]
	\end{equation}
	with the critical intensity $G_c=1.3$. Fig.~\ref{fig3} demonstrates the variation of total intensity $G_s$ throughout the layer of a uniform suspension $(n=1)$ for $V_c=10,\kappa=0.5$, $\omega=0.4$, $I_D=0.25$ and different incidence angles ($\theta_i$). For $0\leq\theta_i\leq 80$, $G_s$ decreases monotonically throughout the suspension and $\theta_i=0$, the critical intensity occurs at mid-height of the suspension. Therefore, cells accumulation in the basic state occurs at the mid-height of the suspension. As $\theta_i$ is further increased, the maximum concentration in the basic state increases and moves nearer to the top of the suspension (see Fig.~\ref{fig3}(b)). In the absence of diffuse irradiation, for $0<\omega<0.7$ the total uniform intensity decreases monotonically with the depth of the suspension. Therefore, the effect of the angle of incidence, $\theta_i$, is the same on concentration in the basic equilibrium state for $0<\omega<0.7$. \par

	\begin{figure*}[!bt]
		\includegraphics{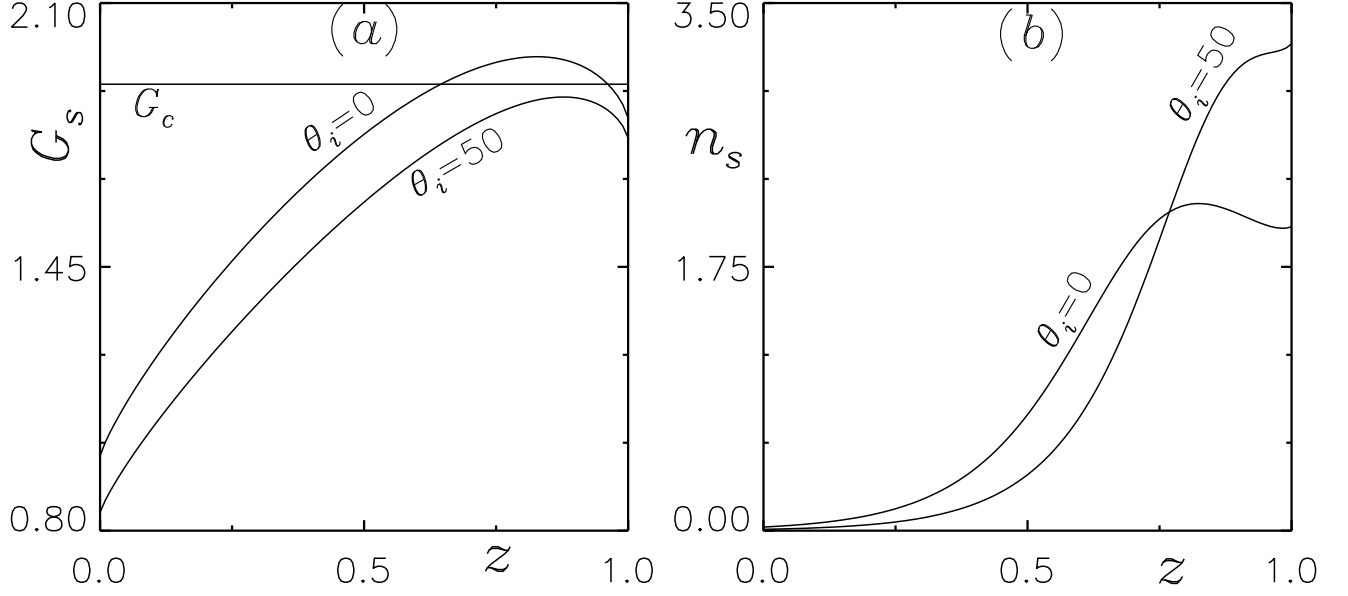}
		\caption{\label{fig4}(a) Variation of total intensity in a uniform suspension for two different values of angle of incidence $\theta_i$ and (b) corresponding base concentration profile for the governing parameter values $S_c=20, V_c=10,k=1, I_D=0.02,\omega=1$ and $I_t=1$ are kept fixed. Here, the critical intensity $G_c=1.9$ is utilized.}
	\end{figure*}

	Now we consider the case of purely scattering suspension ($\omega=1$).  The taxis function
	\begin{equation}\label{35}
		T(G)=0.8\sin\left(\frac{3\pi}{2}\chi(G)\right)-0.1\sin\left(\frac{\pi}{2}\chi(G)\right),\quad \chi(G)=\frac{1}{3.8}G\exp[0.135(3.8-G)]
	\end{equation}
	
	with the critical intensity $G_c=1.9$ is used in this section. The Fig.~\ref{fig4} shows the fluctuation of total intensity and associated basic concentration profile. Here, the total intensity breaks the rule of monotonicity. Therefore, the critical light intensity $G_c$ occurs at two depths $z\approx 0.64$ and $z\approx 0.96$ of the suspension and the bimodal steady state occurs for $\theta_i=0$. On the other hand, for $\theta_i=50$, the basic state occurs at a single location of the suspension (see Fig.~\ref{fig4}).

	\section{Linear stability of the problem}
	Consider a small perturbation of amplitude, $0<\epsilon\leq 1$, to the equilibrium state 
	\begin{widetext}
		
		\begin{align*} 
			[\boldsymbol{u},n,I,<\wp>]=[0,n_s,I_s,<p_s>]+\epsilon [\boldsymbol{u}_1,n_1,I_1,<\boldsymbol{\wp}_1>]+\mathcal{O}(\epsilon^2)=[0,n_s,I_s^c+I_s^d,<\wp_s>]\\
			+\epsilon [\boldsymbol{u}_1,n_1,I_1^c+I_1^d,<\boldsymbol{\wp}_1>]
			+\mathcal{O}(\epsilon^2).  
		\end{align*}
	\end{widetext}
where  $\boldsymbol{u}_1=(u_1,v_1,w_1)$.
	The substitution of perturbed values are made into the Eqs.~(\ref{13})-(\ref{16}) and these equations are made linear about the basic state by collecting $o(\epsilon)$ terms, becomes
	\begin{equation}\label{36}
		\boldsymbol{\nabla}\cdot \boldsymbol{u}_1=0,
	\end{equation}
	
	\begin{equation}\label{37}
		S_{c}^{-1}\left(\frac{\partial \boldsymbol{u_1}}{\partial t}\right)=-\boldsymbol{\nabla} P_{e}+Rn_1\hat{\boldsymbol{z}}-\nabla^{2}\boldsymbol{ u_1},
	\end{equation}
	\begin{equation}\label{38}
		\frac{\partial{n_1}}{\partial{t}}+V_c\boldsymbol{\nabla}\cdot(<\boldsymbol{\wp_s}>n_1+<\boldsymbol{\wp_1}>n_s)+w_1\frac{dn_s}{dz}=\boldsymbol{\nabla}^2n_1.
	\end{equation}
	If $G=G_s+\epsilon G_1+\mathcal{O}(\epsilon^2)=(G_s^c+\epsilon G_1^c)+(G_s^d+\epsilon G_1^d)+\mathcal{O}(\epsilon^2)$, then the total collimated intensity in the basic state is perturbed as $I_t\exp\left(\frac{-\kappa\int_z^1(n_s(z')+\epsilon n_1+\mathcal{O}(\epsilon^2))dz'}{\cos\theta_r}\right)$  and after simplification, we get
	\begin{equation}\label{39}
		G_1^c=I_t\exp\left(\frac{-\int_z^1 \kappa n_s(z')dz'}{\cos\theta_r}\right)\left(\frac{\int_1^z\kappa n_1 dz'}{\cos\theta_r}\right)
	\end{equation}
	and, perturbed basic total diffuse intensity $G_1^d$ can be get by 
	\begin{equation}\label{40}
		G_1^d=\int_0^{4\pi}I_1^d(\boldsymbol{ x},\boldsymbol{ s})d\Omega.
	\end{equation}
	Similarly, for the radiative heat flux $q=q_s+\epsilon q_1++\mathcal{O}(\epsilon^2)=((q_s^c+q_s^d)+\epsilon (q_1^c+q_1^d)+\mathcal{O}(\epsilon^2)$, and we get
	\begin{equation}\label{41}
		\boldsymbol{q}_1^c=-I_t(\cos\theta_r)\exp\left(\frac{-\int_z^1 \kappa n_s(z')dz'}{\cos\theta_r}\right)\left(\frac{\int_1^z\kappa n_1 dz'}{\cos\theta_r}\right)\hat{z}
	\end{equation}
	and
	\begin{equation}\label{42}
		q_1^d=\int_0^{4\pi}I_1^d(\boldsymbol{ x},\boldsymbol{ s})\boldsymbol{ s}d\Omega.
	\end{equation}
	Now the expression 
	\begin{equation*}
		-T(G_s+\epsilon G_1)\frac{\boldsymbol{q}_s+\epsilon\boldsymbol{q}_1+\mathcal{O}(\epsilon^2)}{|\boldsymbol{q}_s+\epsilon\boldsymbol{q}_1+\mathcal{O}(\epsilon^2)|}-T_s\hat{\boldsymbol{z}},
	\end{equation*}
	gives the perturbed swimming orientation by collecting $O(\epsilon)$ terms
	\begin{equation}\label{43}
		<\boldsymbol{\wp_1}>=G_1\frac{dT_s}{dG}\hat{\boldsymbol{z}}-T_s\frac{\boldsymbol{q_1}^H}{\boldsymbol{q_s}},
	\end{equation}
	where $\boldsymbol{q}_1^H$ is the horizontal component of the perturbed radiative heat flux $\boldsymbol{q}_1$.
	Now substituting the value of $<\boldsymbol{\wp_1}>$ from  Eq.~$(\ref{43})$ into Eq.~$(\ref{38})$ and simplifying, we get
	\begin{equation}\label{44}
		\frac{\partial{n_1}}{\partial{t}}+V_c\frac{\partial}{\partial z}\left(T_sn_1+n_sG_1\frac{dT_s}{dG}\right)-V_cn_s\frac{T_s}{q_s}\left(\frac{\partial q_1^x}{\partial x}+\frac{\partial q_1^y}{\partial y}\right)+w_1\frac{dn_s}{dz}=\nabla^2n_1.
	\end{equation}
	Now we eliminate pressure gradient and horizontal component of $u_1$ by taking the double curl Eq.~(\ref{38}) and retaining the z-component of result. Then Eqs.~(\ref{36}~-~\ref{38}) are reduced to two equations for $w_1$ and $n_1$. Now these quantities are decomposed into normal modes such that
	\begin{equation}\label{45}
		w_1=W(z)\exp{(\sigma t+i(lx+my))},\quad n_1=\Phi(z)\exp{(\sigma t+i(lx+my))}.  
	\end{equation}
	The governing equation for perturbed intensity $I_1$ can be written as
	\begin{equation}\label{46}
		\xi\frac{\partial I_1}{\partial x}+\eta\frac{\partial I_1}{\partial y}+\nu\frac{\partial I_1}{\partial z}+\kappa( n_sI_1+n_1I_s)=\frac{\omega\kappa}{4\pi}(n_sG_1+G_sn_1),
	\end{equation}
	with the boundary conditions
	\begin{subequations}
		\begin{equation}\label{47a}
			I_1(x, y, 1, \xi, \eta, \nu) =0,\qquad (\pi/2\leq\theta\leq\pi,0\leq\phi\leq 2\pi), 
		\end{equation}
		\begin{equation}\label{47b}
			I_1(x, y, 0,\xi, \eta, \nu) =0,\qquad (0\leq\theta\leq\pi/2,0\leq\phi\leq 2\pi). 
		\end{equation}
	\end{subequations}
	The $I_1^d$ has the form 
	\begin{equation*}
		I_1^d=\Psi_1^d(z,\xi,\eta,\nu)\exp{(\sigma t+i(lx+my))}. 
	\end{equation*}
	From Eqs.~(\ref{42}) and (\ref{43}), we get
	\begin{equation}\label{48}
		G_1^c=\left[I_t\exp\left(\frac{-\int_z^1 \kappa n_s(z')dz'}{\cos\theta_r}\right)\left(\frac{\int_1^z\kappa n_1 dz'}{\cos\theta_r}\right)\right]\exp{(\sigma t+i(lx+my))}=\mathcal{G}_1^c(z)\exp{(\sigma t+i(lx+my))},
	\end{equation}
	and 
	\begin{equation}\label{49}
		G_1^d=\mathcal{G}_1^d(z)\exp{(\sigma t+i(lx+my))}= \left(\int_0^{4\pi}\Psi_1^d(z,\xi,\eta,\nu)d\Omega\right)\exp{(\sigma t+i(lx+my))},
	\end{equation}

	where $\mathcal{G}_1(z)=\mathcal{G}_1^c(z)+\mathcal{G}_1^d(z)$ is the perturbed total intensity.
	
	Now $\Psi_1^d$ satisfies
	\begin{equation}\label{50}
		\frac{d\Psi_1^d}{dz}+\frac{(i(l\xi+m\eta)+\kappa n_s)}{\nu}\Psi_1^d=\frac{\omega\kappa}{4\pi\nu}(n_s\mathcal{G}_1+G_s\Phi)-\frac{\kappa}{\nu}I_s\Phi,
	\end{equation}
	subject to the boundary conditions
	\begin{subequations}
		\begin{equation}\label{51a}
			\Psi_1^d( 1, \xi, \eta, \nu) =0,\qquad (\pi/2\leq\theta\leq\pi,0\leq\phi\leq 2\pi), 
		\end{equation}
		\begin{equation}\label{51b}
			\Psi_1^d( 0,\xi, \eta, \nu) =0,\qquad (0\leq\theta\leq\pi/2,0\leq\phi\leq 2\pi). 
		\end{equation}
	\end{subequations}
	Similarly from Eq.~(\ref{43}), we have
	\begin{equation*}
		q_1^H=[q_1^x,q_1^y]=[P(z),Q(z)]\exp{[\sigma t+i(lx+my)]},
	\end{equation*}
	where
	\begin{equation*}
		P(z)=\int_0^{4\pi}\Psi_1^d(z,\xi,\eta,\nu)\xi d\Omega,\quad Q(z)=\int_0^{4\pi}\Psi_1^d(z,\xi,\eta,\nu)\eta d\Omega.
	\end{equation*}
	The linear stability equations become
	\begin{equation}\label{52}
		\left(\sigma S_c^{-1}+k^2-\frac{d^2}{dz^2}\right)\left( \frac{d^2}{dz^2}-k^2\right)W=Rk^2\Phi,
	\end{equation}
	\begin{equation}\label{53}
		\left(\sigma+k^2-\frac{d^2}{dz^2}\right)\Phi+V_c\frac{d}{dz}\left(T_s\Phi+n_s\mathcal{G}_1\frac{dT_s}{dG}\right)-i\frac{V_cn_sT_s}{q_s}(lP+mQ)=-\frac{dn_s}{dz}W,
	\end{equation}
	subject to the boundary conditions
	\begin{equation}\label{54}
		W=\frac{d^2W}{dz^2}=\frac{d\Phi}{dz}-V_cT_s\Phi-n_sV_C\mathcal{G}_1\frac{dT_s}{dG}=0,\quad at\quad z=0,1.
	\end{equation}
	
	Here, $k=\sqrt{(l^2+m^2)}$ is the overall non-dimensional wavenumber. Eqs.~(\ref{52})-(\ref{53}) form an eigen value problem for $\sigma$ as a function of the dimensionless parameters $V_c,\kappa,\sigma,l,m,R$.
	
	Eq.~(\ref{53}) becomes 
	\begin{equation}\label{55}
		\Gamma_0(z)+\Gamma_1(z)\int_1^z\Theta dz+(\sigma+k^2+\Gamma_2(z))\Phi+V_cT_sD\Phi-D^2\Phi=-Dn_sW, 
	\end{equation}
	where
	\begin{subequations}
		\begin{equation}\label{56a}
			\Gamma_0(z)=V_cD\left(n_s\mathcal{G}_1^d\frac{dT_s}{dG}\right)-i\frac{V_cn_sT_s}{q_s}(lP+mQ),
		\end{equation}
		\begin{equation}\label{56b}
			\Gamma_1(z)=\kappa V_cD\left(n_sG_s^c\frac{dT_s}{dG}\right)
		\end{equation}
		\begin{equation}\label{56c}
			\Gamma_2(z)=2\kappa V_c n_s G_s^c\frac{dT_s}{dG}+V_c\frac{dT_s}{dG}DG_s^d,
		\end{equation}
	\end{subequations}

	Now introducing a new variable
	\begin{equation}\label{57}
		\Theta(z)=\int_1^z\Phi(z')dz',
	\end{equation}
	the linear stability equations become
	\begin{equation}\label{58}
		\left(\sigma S_c^{-1}+k^2-D^2\right)\left( D^2-k^2\right)W=Rk^2D\Theta,
	\end{equation}
	\begin{equation}\label{59}
		\Gamma_0(z)+\Gamma_1(z)\Theta+(\sigma+k^2+\Gamma_2(z))D\Theta+V_cT_sD^2\Theta-D^3\Theta=-Dn_sW. 
	\end{equation}
	The boundary conditions become,
	
	\begin{equation}\label{60}
		W=D^2W=D^2\Theta-\Gamma_2(z)D\Theta-V_cT_s\frac{dT_s}{dG}\mathcal{G}_1=0,\quad at\quad z=0,1.
	\end{equation}
and the additional boundary condition is given by
	\begin{equation}\label{61}
	\Theta(z)=0,\quad at\quad z=1.
	\end{equation}
	
which is augmented by the Eq.~\ref{57}.	
	
	\section{SOLUTION PROCEDURE}
To solve Eqs.~(\ref{58}) and (\ref{59}) with appropriate boundary conditions, a finite-difference scheme based on Newton-Raphson-Kantorovich iterations is used, with fourth-order accuracy. This scheme yields neutral stability curves for a fixed set of other parameters. Initial values for $S_c, V_c, \kappa, \omega, k, \theta_i$, and $I_D$ are given, and W and $\Theta$ are estimated either from previous results or by imposing sinusoidal variation. Once a solution is obtained, it can be used as an initial guess for nearby parameter values. For each set of parameter ranges, there are infinitely many branches of the neutral curve $R^{(n)}(k)$, each representing a different solution of the linear stability problem. The most interesting branch is the one where R has its minimum value, $R_c$, and the most unstable solution is the pair $(k_c, R_c)$, which gives the wavelength of the initial disturbance, $\lambda_c = 2\pi/k_c$. Bioconvective solutions consist of convection cells stacked vertically in the suspension, with a mode n solution having n such cells. Often, the most unstable solution is found on the $R^{(1)}(k)$ branch of the neutral curve, which is mode 1. A neutral curve is where Re$(\sigma) = 0$, and if Im$(\sigma) = 0$ on such a curve, the bioconvective solution is stationary. If Im$(\sigma)\neq 0$, then oscillatory solutions exist, and if the most unstable solution remains on the oscillatory branch, it's overstable. Oscillatory solutions arise when there is a competition between stabilizing and destabilizing processes. When oscillatory solutions occur, a single oscillatory branch of the neutral curve meets the stationary branch at $k = k_0$ and exists for $k \leq k_0$.
		
	\section{NUMERICAL RESULTS}
	In this article, the effect of angle of incidence $\theta_i$ is investigated where other governing parameters $S_c,I_t,V_c,\kappa,I_D$, and $\omega$ are kept constant. In the presence of various parameter values, it is hard to get a complete picture for whole domain.
	So, a discrete set of fixed parameter values are consider. The value of schmidt number $S_c=20$ and total direct flux (collimated )$I_t=1$ are kept constants throughout the study. Other parameters like scattering albedo, extinction coefficient, and cell swimming speed are taken as $\omega\in [0: 1]$, $\kappa= 0.5, 1.0$, and $V_c=10,15,20$. Here, we differentiate the study into two cases based on the effectiveness of self-shading and scattering.
	
	\subsection{WHEN SELF SHADING IS EFFECTIVE}
	
	To study the effect of incidence angle, $\theta_i$, on bioconvective instability, first we consider the case when self-shading is effective. In this section scattering is considered weak by selecting the lower value of scattering albedo $\omega$ and here we discuss two cases when $\kappa=$ 0.5 and $\kappa=$1. The critical intensity $G_c=1.3$ is utilized throughout the study.
	\begin{figure*}[b]
		\includegraphics{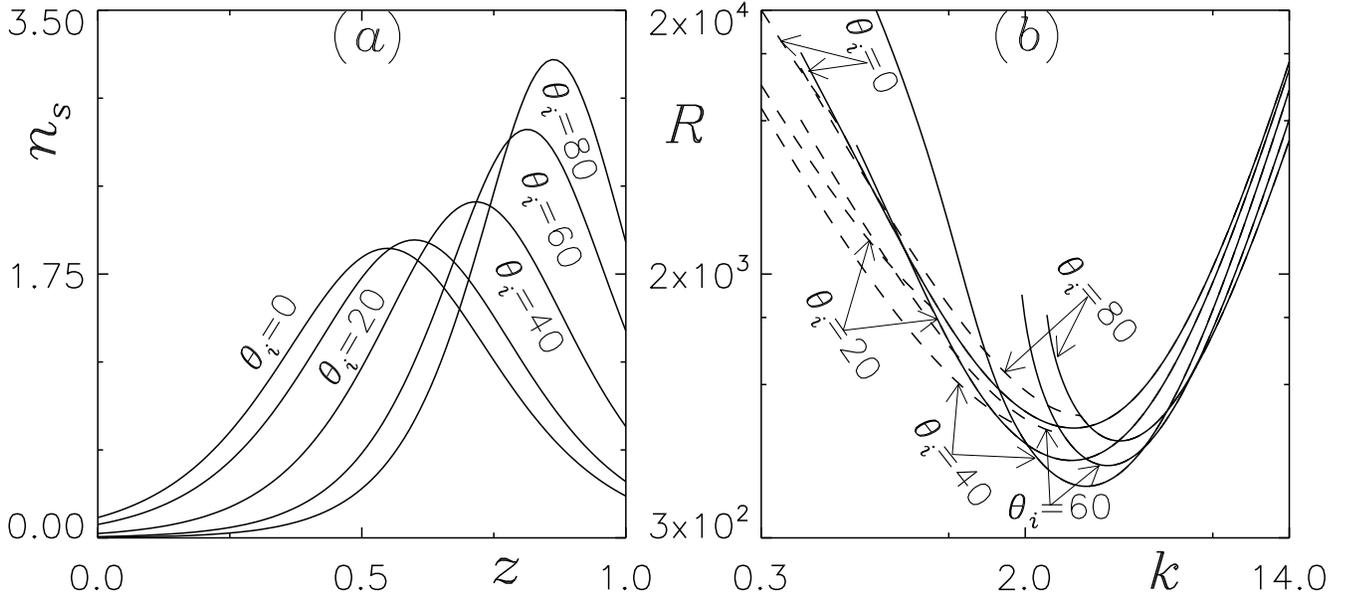}
		\caption{\label{fig5}(a) Basic concentration profile, (b) corresponding marginal stability curves for different values of angle of incidence $\theta_i$. Here, the governing parameter values $S_c=20,V_c=15,k=0.5, I_D=0.26$, and $\omega=0.4$ are kept fixed.}
	\end{figure*}
	
	\subsubsection{$V_c=15$}
	(\romannumeral 1) When extinction coefficient $\kappa=0.5$\\
	The basic cell concentration profile and corresponding neutral curves are shown in fig.~\ref{fig5} for different values of incidence angle $\theta_i$, where the other parameters $v_c=15,k=0.5,I_D=0.26$, and $\omega=0.4$ are kept fixed. When $\theta_i=0$, the maximum cell concentration occurs around at domain's mid-hight. As the value of $\theta_i$ increases to 20, the maximum concentration shifts at $z\approx 0.6$, and the unstable region's width increases, supporting the convective fluid motion. As a result, the lower critical Rayleigh number $R_c$ occurs. The basic state occurs at $z\approx 0.72$ for $\theta_i=40$ and the width of the unstable region increases. Here, due to an increment in an unstable region, the critical Rayleigh number decreases compared to the case of $\theta_i=20$. When $\theta_i=60$, the cells accumulate at $z\approx$ 0.8. Here, due to an increment in an angle of incidence $\theta_i$, the width of the unstable region and steepness of the maximum concentration increases, which supports the convection, but on the other hand, resistance due to positive phototaxis also increases, which inhibit the convective fluid motion. In this case, the latter effect dominates the former, and the higher Rayleigh number occurs. As $\theta_i$ increases, the location of the maximum concentration shifts towards the domain's top, and similar effects are seen on the critical Rayleigh number. In all the cases, one oscillatory branch bifurcates from the stationary branch of the marginal (neutral) stability curves and exists throughout less than bifurcation points. However, the most unstable solution occurs on the stationary branch of the marginal stability curves. Therefore, perturbation to the basic state remains stationary throughout this section.   
	
		\begin{figure*}[b]
		\includegraphics{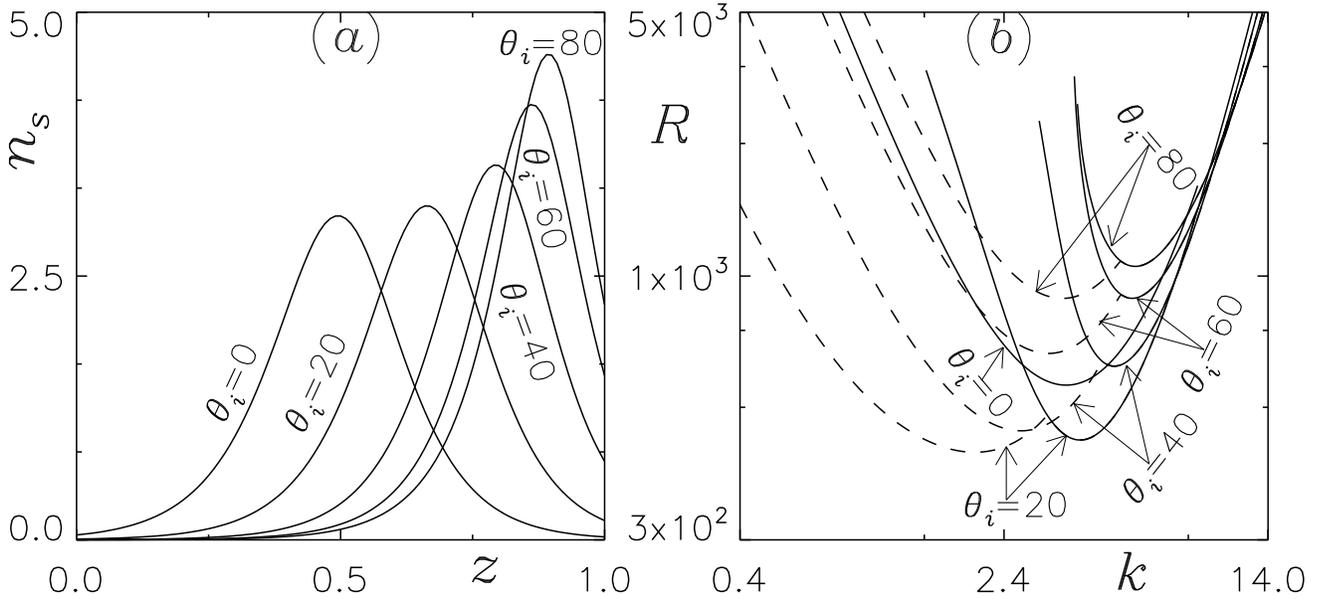}
		\caption{\label{fig6}(a)Basic concentration profile, (b) corresponding marginal stability curves for different values of angle of incidence $\theta_i$. Here, the governing parameter values $S_c=20,V_c=15,k=1, I_D=0.5$, and $\omega=0.4$ are kept fixed.}
	\end{figure*}
	
	(\romannumeral 2) When extinction coefficient $\kappa=1$\\
	Fig.~\ref{fig6} shows the effect of increment in an angle of incidence $\theta_i$ on the basic state and corresponding marginal stability curves. Here, the governing parameters $V_c=15,\kappa=1,\omega=0.4$ and $I_D=0.5$ are kept fixed.
	At $\theta_i=0$, the location of the maximum basic concentration is near the domain's mid-hight, and the most unstable bioconvective solution remains in the stationary branch leading the solution to be stationary (non-oscillatory). As $\theta_i$ increases, the location of the maximum basic concentration
	shifts toward the domain's top. The location of the maximum basic concentration occurs at around $z\approx 0.66$ for $\theta=20$. Here. An oscillatory branch bifurcates from the stationary branch of the marginal stability curve at $k\approx$3.18, which exists throughout $k\leq 3.18$. In this case, the most unstable solution occurs on the oscillatory branch of the marginal curve. Therefore, the perturbation to the basic state becomes overstable for $\theta_i=20$. For $\theta_i=40$, the cells accumulate at $z\approx 0.79$ in the basic state, and an oscillatory branch splits from the stationary branch at $k\approx 4.59$, and it remains throughout $k\leq 4.59$. Here, the most unstable solution is overstable, and the perturbation to the basic state remains oscillatory. When $\theta_i$ is increased to 60, the location of the cell accumulation in the basic state shift at $z\approx$ 0.86, and here the oscillatory branch of the marginal curve bifurcate from the stationary branch at $k\approx5.21$ which also remains all over $k\leq 5.21$. The oscillatory branch has the most unstable solution, and perturbation remains oscillatory here. As $\theta_i$ is further increased to 80, the location of the maximum cell concentration shifts nearer to the domain's top. Here the neutral curve shows the same behavior as an oscillatory branch of the neutral curve bifurcates from the stationary branch, and the most unstable solution occurs on the oscillatory branch. So perturbation to the basic state remains oscillatory here also.
	The quantitative results of this section are shown in Table~\ref{tab1}.
		
	\begin{table}[h]
		\caption{\label{tab1}The quantitative values of bioconvective solutions for discrete values of angle of incidence $\theta_i$ for $V_c=15$ are shown in the table, where other parameters are kept fixed.}
		\begin{ruledtabular}
			\begin{tabular}{cccccccc}
				$V_c$ & $\kappa$ & $\omega$ & $I_D$ & $\theta_i$ & $\lambda_c$ & $R_c$ & $Im(\sigma)$ \\
				\hline
				15 & 0.5 & 0.4 & 0.26 & 0\footnotemark[1]  & 2.74 & 719.02 & 0 \\
				15 & 0.5 & 0.4 & 0.26 & 20\footnotemark[1] & 2.18 & 555.59 & 0 \\
				15 & 0.5 & 0.4 & 0.26 & 40\footnotemark[1] & 1.95 & 452.24 & 0 \\
				15 & 0.5 & 0.4 & 0.26 & 60\footnotemark[1] & 1.68 & 533.648 & 0 \\
				15 & 0.5 & 0.4 & 0.26 & 80\footnotemark[1] & 1.54 & 648.88 & 0 \\
				15 & 1 & 0.4 & 0.5 & 0  & 1.76 & 684.66  & 0 \\
				15 & 1 & 0.4 & 0.5 & 20  & 3.23\footnotemark[2] & 478.61\footnotemark[2] & 11.7 \\
				15 & 1 & 0.4 & 0.5 & 40 & 2.35\footnotemark[2] & 536.69\footnotemark[2] & 18.7 \\
				15 & 1 & 0.4 & 0.5 & 60 & 1.98\footnotemark[2] & 810.67\footnotemark[2] & 23.39 \\
				15 & 1 & 0.4 & 0.5 & 80 & 1.78\footnotemark[2] & 1087.01\footnotemark[2] & 23.11 \\
			\end{tabular}
		\end{ruledtabular}
		\footnotetext[1]{A result indicates that the $R^{(1)}(k)$ branch of the neutral curve is oscillatory.}
		\footnotetext[2]{A result indicates that a smaller solution occurs on the oscillatory branch.}
	\end{table}
	
		\begin{figure*}[b]
		\includegraphics{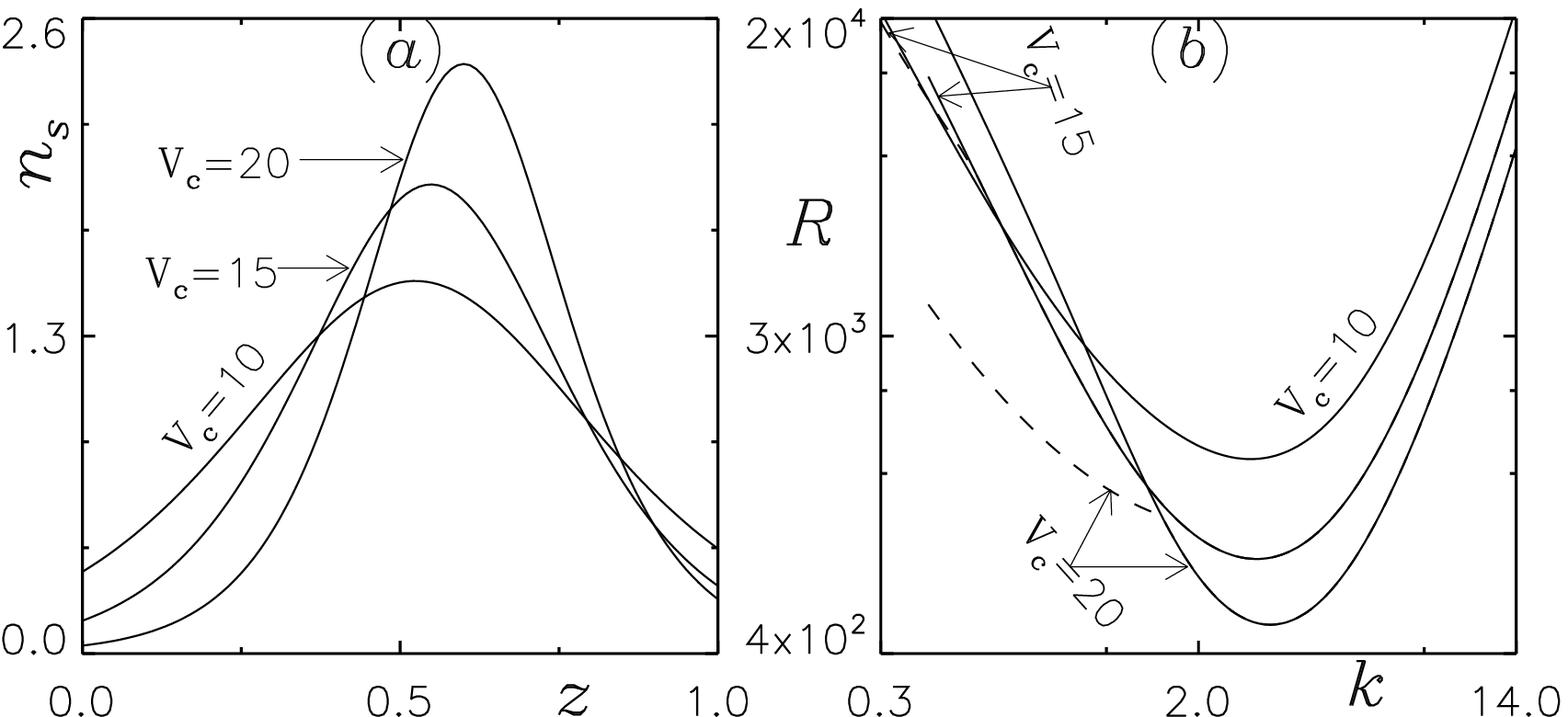}
		\caption{\label{fig7}(a) Basic concentration profiles and (b) corresponding marginal stability curves for variation in swimming speed $V_c$. Here, other governing parameter values $S_c=20,I_D=0.26,k=0.5,\theta_i=0$, and $\omega=0.4$ are kept fixed.}
	\end{figure*}

	\subsubsection{$V_c=10$ and $V_c=20$}
	We have also investigated the effect of the angle of incidence on the bio-convective instability for $V_c=10$ and $V_c=20$. Table~\ref{tab2} provides an overview of the numerical data for the critical Rayleigh number and wave number for $V_c=10$ and $V_c=20$.

	\begin{table}[h]
		\caption{\label{tab2}The quantitative values of bioconvective solutions for discrete values of angle of incidence $\theta_i$ for $V_c=10$ and $V_c=20$ are shown in the table, where other parameters are kept fixed.}
		\begin{ruledtabular}
			\begin{tabular}{cccccccc}
				$V_c$ & $\kappa$ & $\omega$ & $I_D$ & $\theta_i$ & $\lambda_c$ & $R_c$ & $Im(\sigma)$ \\
				\hline
				10 & 0.5 & 0.4 & 0.25 & 0  & 2.51 & 942.81 & 0 \\
				10 & 0.5 & 0.4 & 0.25 & 20 & 2.63 & 775.87 & 0 \\
				10 & 0.5 & 0.4 & 0.25 & 40& 2.63 & 535.98 & 0 \\
				10 & 0.5 & 0.4 & 0.25 & 60 & 2.35 & 451.21 & 0 \\
				10 & 0.5 & 0.4 & 0.25 & 80 & 2.13 & 457.42 & 0 \\
				10 & 1 & 0.4 & 0.5 & 0  & 1.92 & 860.48  & 0 \\
				10 & 1 & 0.4 & 0.5 & 20  & 1.95 & 637.64 & 11.7 \\
				10 & 1 & 0.4 & 0.5 & 40\footnotemark[1] & 1.76 & 521.08 & 18.7 \\
				10 & 1 & 0.4 & 0.5 & 60\footnotemark[1] & 1.56 & 603.81 & 23.39 \\
				10 & 1 & 0.4 & 0.5 & 80\footnotemark[1] & 1.52 & 683.46 & 23.11 \\
				20 & 0.5 & 0.4 & 0.25 & 0\footnotemark[1]  & 1.95 & 666.93 & 0 \\
				20 & 0.5 & 0.4 & 0.25 & 20\footnotemark[1] & 1.94 & 461.48 & 0 \\
				20 & 0.5 & 0.4 & 0.25 & 40\footnotemark[1] & 1.58 & 523.57 & 0 \\
				20 & 0.5 & 0.4 & 0.25 & 60 & 2.01\footnotemark[2] & 758.31\footnotemark[2] & 13.95 \\
				20 & 0.5 & 0.4 & 0.25 & 80 & 1.79\footnotemark[2] & 1004.45\footnotemark[2] & 16.49 \\
				20 & 1 & 0.4 & 0.5 & 0\footnotemark[1]  & 1.95 & 860.48  & 0 \\
				20 & 1 & 0.4 & 0.5 & 20  & 2.9\footnotemark[2] & 637.64\footnotemark[2] & 19.58 \\
				20 & 1 & 0.4 & 0.5 & 40 & 2.09\footnotemark[2] & 521.08\footnotemark[2] & 32.34 \\
				20 & 1 & 0.4 & 0.5 & 60 & 1.7\footnotemark[2] & 603.81\footnotemark[2] & 41.13 \\
				20 & 1 & 0.4 & 0.5 & 80 & 1.54\footnotemark[2] & 683.46\footnotemark[2] & 41.62 \\	
			\end{tabular}
		\end{ruledtabular}
		\footnotetext[1]{A result indicates that the $R^{(1)}(k)$ branch of the neutral curve is oscillatory.}
		\footnotetext[2]{A result indicates that a smaller solution occurs on the oscillatory branch.}
	\end{table}

	\subsubsection{Effect of swimming speed}
\begin{figure*}[!bt]
		\includegraphics{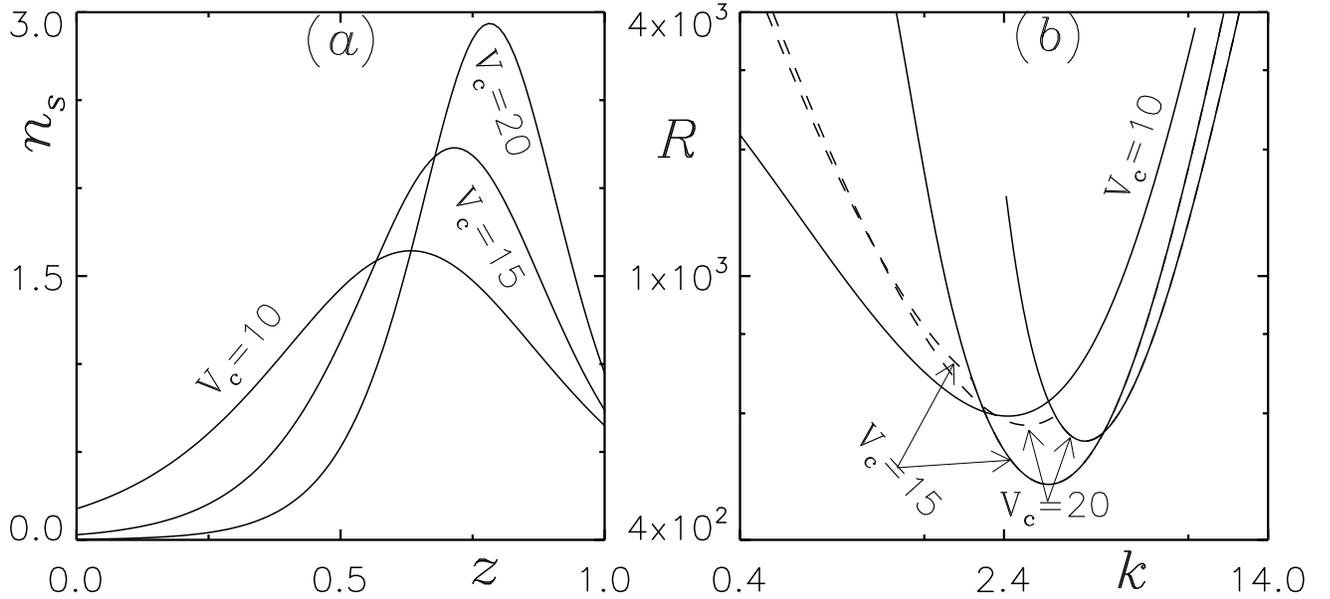}
		\caption{\label{fig8}(a) Basic concentration profiles and (b) corresponding marginal stability curves for variation in swimming speed $V_c$. Here, other governing parameter values $S_c=20,I_D=0.26,\theta_i=40,k=0.5$, and $\omega=0.4$ are kept fixed.}
	\end{figure*}

\begin{figure*}[!bt]
		\includegraphics{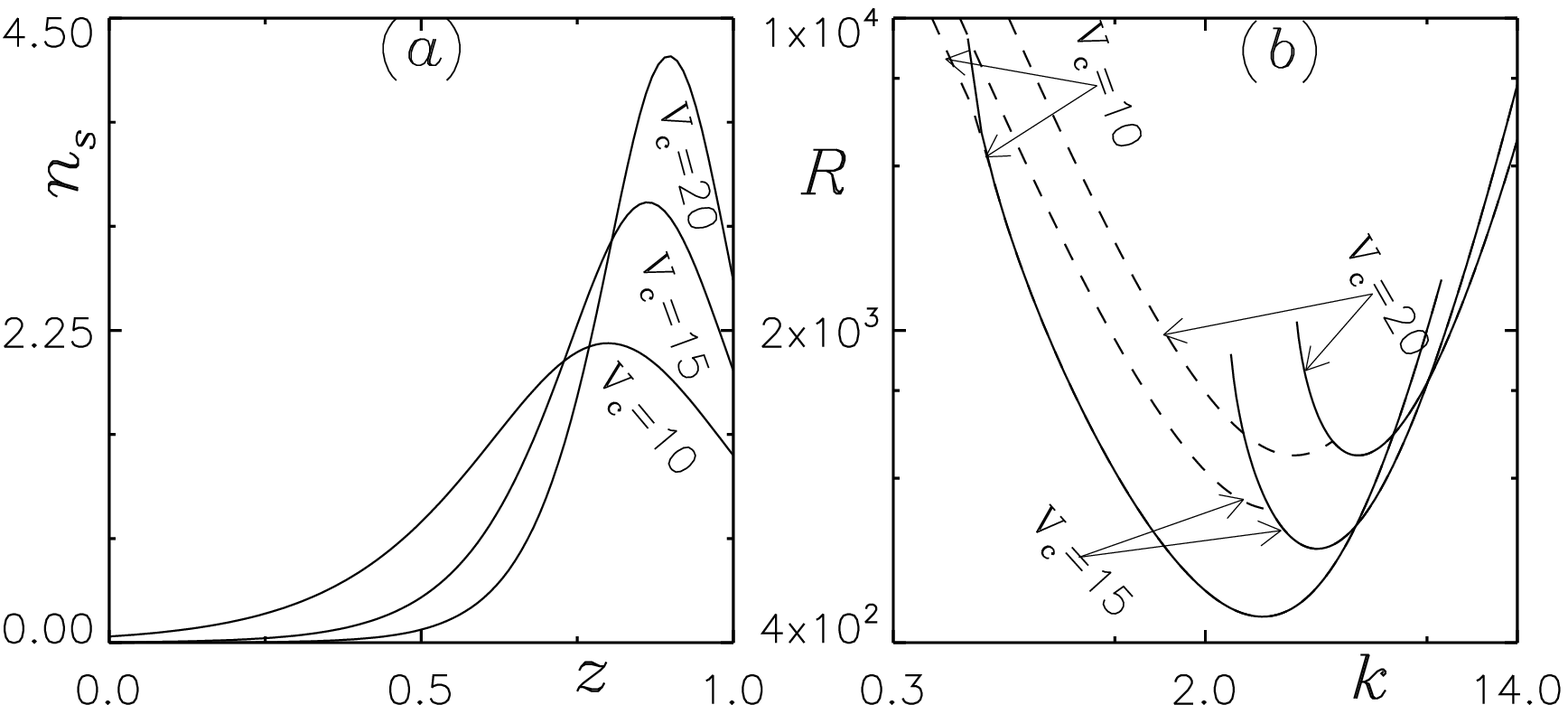}
		\caption{\label{fig9}(a) Basic concentration profiles and (b) corresponding marginal stability curves for variation in swimming speed $V_c$. Here, other governing parameter values $S_c=20,I_D=0.26,k=0.5,\theta_i=80$, and $\omega=0.4$ are kept fixed.}
	\end{figure*}

		To study the effect of cell swimming speed $V_c$ on critical values of the Rayleigh number and wavelength, we vary the value of swimming speed from 10 to 20. The basic profile of cell concentration and associated marginal stability curves are shown in  Figs.~\ref{fig7} to \ref{fig9}, for three cases $\theta_i=0,40,80$, where the other governing parameters $\kappa=0.5,\omega=0.4$ and $I_D=0.26$ are kept fixed.
	The first case shows that as the swimming speed $V_c$ increases, the concentration at top of suspension becomes steeper, supporting the convective fluid motion. However, at higher swimming speeds, the positive phototaxis causes cells in the plum to face higher resistance, which inhibits the convection. Therefore, the Rayleigh number decreases with an increase in $V_c$. In the second case, the maximum concentration occurs at a higher location and increases as $V_c$ increases, but again, positive phototaxis hinders the convective fluid motion, resulting in a higher critical Rayleigh number for higher swimming speeds. An oscillatory branch bifurcates from the stationary branch for $V_c=$ 15 and 20, but the most unstable solution occurs on the stationary branch. In the last case, the concentration becomes steeper as $V_c$ increases, supporting the convective fluid motion, but the region of positive phototaxis increases as well, which opposes the convection. Nonetheless, the former effect dominates, resulting in a lower critical Rayleigh number. An oscillatory branch is seen here too, but the most unstable solution occurs on the stationary branch, and perturbation to the basic state is stationary in all cases.

Here, the impact of cell swimming speed on a bio-convective solution has been also examined for a value of $\kappa=$ 1. The effect of swimming speed is found to be different from that observed for $\kappa$ equal to 0.5. Specifically, the critical Rayleigh numbers are observed to increase for higher swimming speeds, which contrasts with the results for $\kappa=$0.5 and $\theta_i=$ 40. Additionally, the bio-convective solution is found to be frequently overstable for higher values of swimming speed $V_c$. Table~\ref{tab3} presents the numerical outcomes for the bioconvective solutions.
	
		\begin{table}
		\caption{\label{tab3}The quantitative values of bioconvective solutions showing the effect of variation in cell swimming speed are shown in the table, where other parameters are kept fixed.}
		\begin{ruledtabular}
			\begin{tabular}{ c c c c c c c c}
				$\theta_i$ & $I_D$  & $\kappa$ & $\omega$ & $V_c$ & $\lambda_c$ & $R_c$ & $I_m(\sigma)$ \\
				\hline
				
				0 & 0.26 & 0.5 & 0.4 & 10 & 2.23 & 1335.29  & 0 \\
				0 & 0.26 & 0.5 & 0.4 & 15\footnotemark[1] & 2.14 & 719.02 & 0 \\
				0 & 0.26 & 0.5 & 0.4 & 20\footnotemark[1] & 1.98 & 478.47 & 0 \\
				40 & 0.26 & 0.5 & 0.4 & 10 & 2.61 & 619.83 & 0\\			
				40 & 0.26 & 0.5 & 0.4 & 15\footnotemark[1] & 1.95 & 452.24 & 0 \\
				40 & 0.26 & 0.5 & 0.4 & 20\footnotemark[1] & 1.54 & 552.18 & 0 \\
				80 & 0.26 & 0.5 & 0.4 & 10\footnotemark[1] & 2.14 & 457.35 & 0\\
				80 & 0.26 & 0.5 & 0.4 & 15\footnotemark[1] & 1.54 & 648.88 & 0\\
				80 & 0.26 & 0.5 & 0.4 & 20 & 1.79\footnotemark[2] & 1050.19\footnotemark[2] & 16.24 \\
				0 & 0.48 & 1 & 0.4 & 10\footnotemark[1] & 1.98 & 592.73  & 0\\
				0 & 0.48 & 1 & 0.4 & 15 & 2.97\footnotemark[2] & 459.75\footnotemark[2] & 12.57 \\
				0 & 0.48 & 1 & 0.4 & 20 & 2.69\footnotemark[2] & 423.34\footnotemark[2] & 21.52 \\
				40 & 0.48 & 1 & 0.4 & 10\footnotemark[1] & 1.70 & 533.14 & 0 \\
				40 & 0.48 & 1 & 0.4 & 15 & 2.21\footnotemark[2] & 618.01\footnotemark[2] & 20.06 \\
				40 & 0.48 & 1 & 0.4 & 20 & 1.91\footnotemark[2] & 820.63\footnotemark[2] & 35.22 \\
				80 & 0.48 & 1 & 0.4 & 10\footnotemark[1] & 1.98 & 915.67 & 0\\
				80 & 0.48 & 1 & 0.4 & 15 & 1.76\footnotemark[2] & 1190.16\footnotemark[2] & 21.13\\
				80 & 0.48 & 1 & 0.4 & 20 & 1.5\footnotemark[2] & 1854.45\footnotemark[2] & 38.84 \\
				
			\end{tabular}
		\end{ruledtabular}
		\footnotetext[1]{A result indicates that the $R^{(1)}(k)$ branch of the neutral curve is oscillatory.}
		\footnotetext[2]{A result indicates that a smaller solution occurs on the oscillatory branch.}
	\end{table}

	\subsection{WHEN SCATTERING IS EFFECTIVE}
	
	\begin{figure*}[!bt]
		\includegraphics{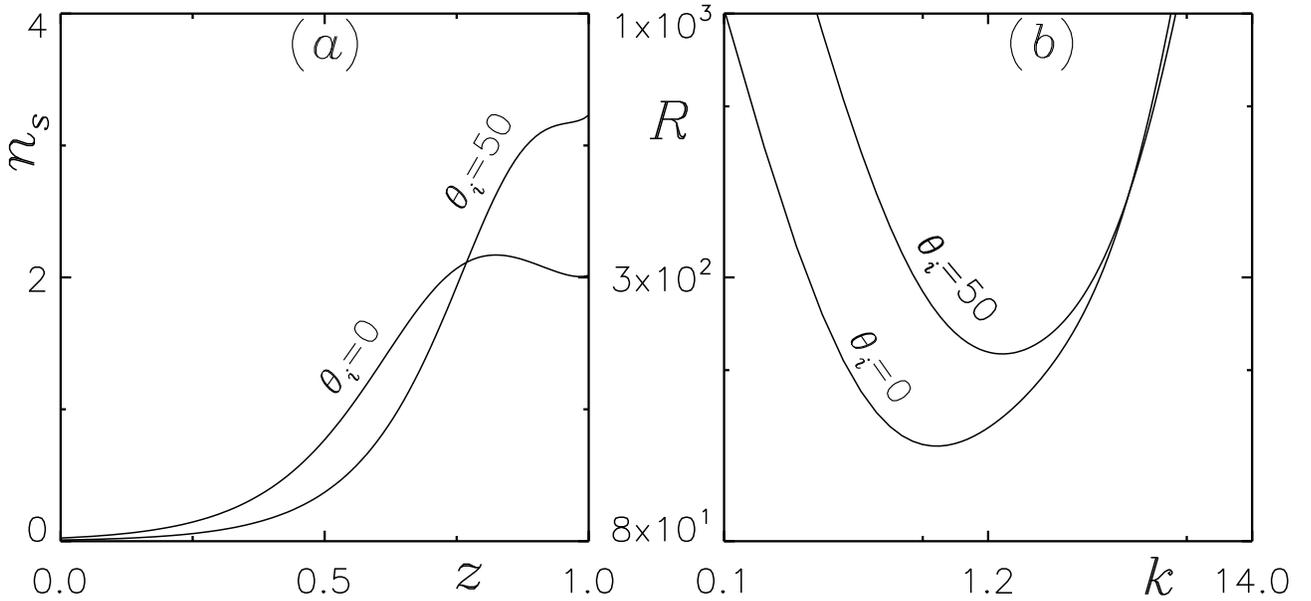}
		\caption{\label{fig10}(a) Basic concentration profiles and (b) corresponding marginal stability curves for $\theta=0$ and 50. Here, the suspension is assumed to be almost purely scattering ($\omega=1$), and other governing parameter values $S_c=20,k=1, I_D=0.02$, and $V_c=10$ are kept fixed.}
	\end{figure*}

Now we investigate the influence of the angle of incidence ($\theta_i$) on critical wavenumber and Rayleigh number at bioconvective instability for purely scattering suspensions ($\omega=1$). We assume $\kappa=1$, $\omega=1$, and $I_D=0.02$ and consider three different cases with $V_c=10, 15,$ and $20$. To study the effect of $\theta_i$, we vary it from 0 to 40 and observe its impact on the basic state and the corresponding critical values ($R_c$ and $\lambda_c$). The phototaxis function, similar to Eq.(\ref{37}), is employed with the critical intensity $G_c=1.9$ in each case. Self-shading is negligible as we consider $\omega=1$.
	\begin{figure*}[!bt]
		\includegraphics{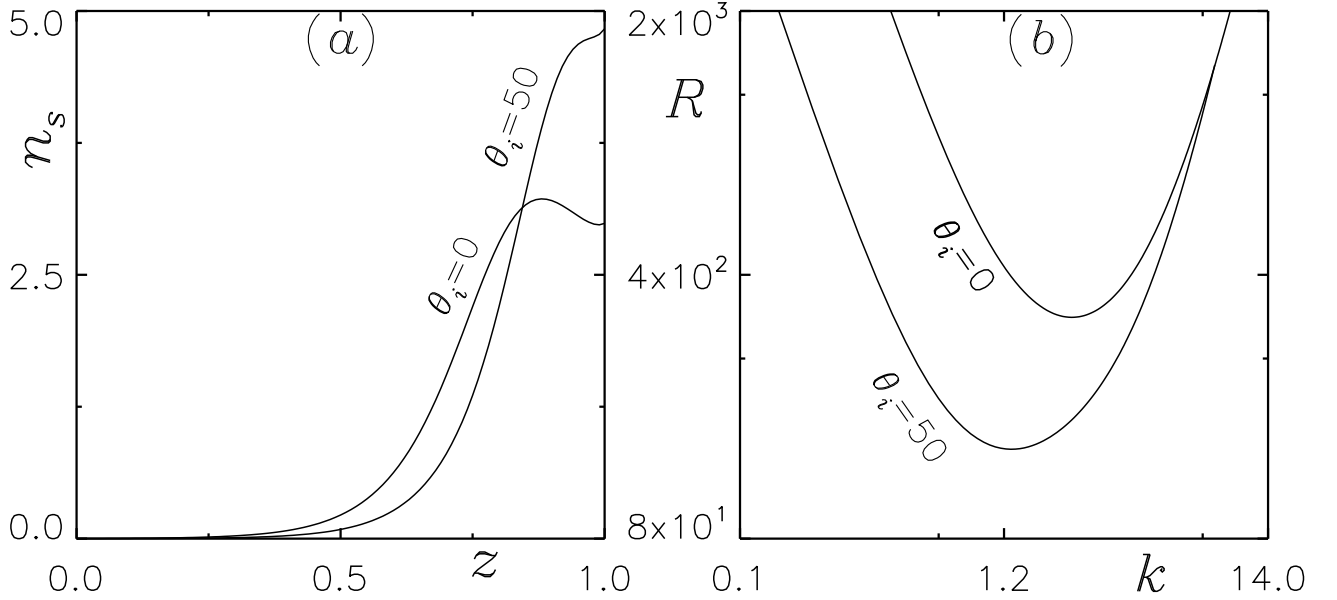}
		\caption{\label{fig11}(a) Basic concentration profiles and (b) corresponding marginal stability curves for $\theta=0$ and 50. Here, the suspension is assumed to be almost purely scattering ($\omega=1$), and other governing parameter values $S_c=20,k=1, I_D=0.02$, and $V_c=15$ are kept fixed.}
	\end{figure*}

\begin{figure*}[!bt]
		\includegraphics{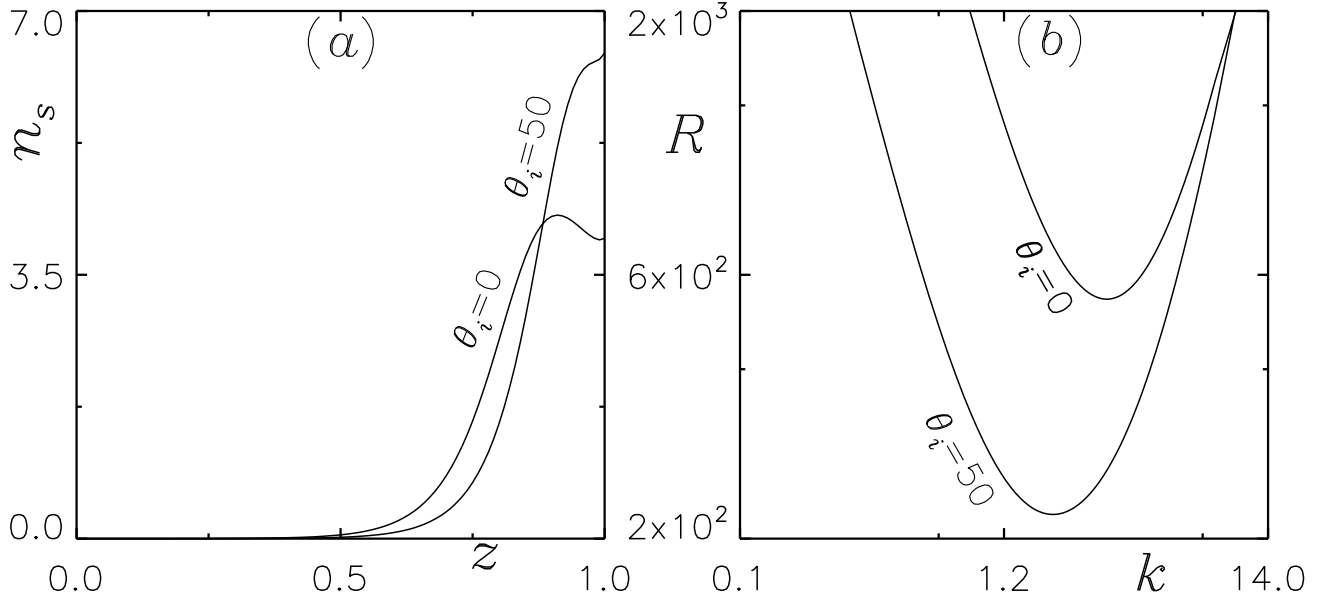}
		\caption{\label{fig12}(a) Basic concentration profiles and (b) corresponding marginal stability curves for $\theta=0$ and 50. Here, the suspension is assumed to be almost purely scattering ($\omega=1$), and other governing parameter values $S_c=20,k=1, I_D=0.02$, and $V_c=20$ are kept fixed.}
	\end{figure*}

Figs.~\ref{fig10}, \ref{fig11}, and \ref{fig12} depict the cell concentration in the basic state and marginal stability curves for three cases with $V_c=10$, 15, and 20 at $\theta_i=0$ and 50, where $S_c=20$, $\omega=1$, $\kappa=1$, and $I_D=0.02$. At $\theta_i=0$, the basic bimodal state is observed in all cases at two locations in the medium. Positive phototaxis occurs below the lower and above the upper locations, while negative phototaxis occurs in between these two locations, causing cell accumulation at two locations closer to the top of the suspension. As $\theta_i$ is increased to 50, all or some algae cells between the two locations swim upward due to positive phototaxis resulting from dim light availability. Thus, cell accumulation occurs in the basic state at the top of the suspension, and the bimodal steady state becomes unimodal. The unstable region width increases with an increase in $\theta_i$ in all cases, and the steepness of higher concentration in the basic state also increases, reinforces the convective fluid motion. However, an increase in $V_c$ from 10 to 20 leads to cells in the residing plum facing higher resistance due to positive phototaxis, which opposes convection. The first effect dominates the latter for higher swimming speeds, resulting in a lower critical Rayleigh number for $\theta_i=50$ compared to $\theta_i=0$ for $V_c=15$ and $V_c=20$. 
		
\section{Conclusion}
	This article explores the impact of a rigid top surface with both diffuse and collimated oblique irradiation on the onset of phototactic bioconvection in a suspension of isotropic scattering phototactic algae. Here, we perturbed the basic equilibrium state and check the linear stability of the same suspension by using the linear perturbation theory similar to the previous published works.

	Isotropic scattering has a significant impact on total uniform intensity as well basic state in the algal suspension. In the case of purely scattering suspension, the variation of total intensity is not monotonic throughout the suspension depth. First, it increases and then decreases. So, the critical intensity occurs at the two locations in the suspension's domain. As a result, the bimodal steady state is observed for purely scattering suspension, which converts into a unimodal steady state as the angle of incidence increases.
	
	The linear stability analysis shows the both types of the nature of disturbance (stationary and oscillatory). Oscillatory behavior of the solution occurs due to conflict between the processes of stabilizing and destabilizing in the suspension. We have also observed that the bioconvective solutions transit from stationary to oscillatory and vice versa as the angle of incidence increases with the fixed governing parameter values. The oscillatory solutions are also observed as swimming speed increases.
	
The presence of a rigid top surface causes the fluid flow near the surface to be inhibited. Consequently, a larger Rayleigh number required for the convective fluid motion in the suspension. This finding is supported by numerical data, which are shown in Tables \ref{tab1}, \ref{tab2}, and \ref{tab3}. The critical Rayleigh number rises as the diffuse irradiation intensity and angle of incidence both also increases.
As a result, when a suspension has rigid top surface  which illuminates by both diffuse and collimated oblique irradiation, the suspension becomes more stable. 
	
It is important to compare the theoretical predictions with the quantitative experimental findings for the suspension of purely  phototactic algae but till the present date there aren't any such statistics accessible. Therefore, to search the purely phototactic algal species are very important but existing knowledge about species are shown that the mostly algal species which shows the phototactic nature is gravitactic and gyrotactic also. This model can be used to handle the other different problems such as anisotropic scattering by changing the scattering phase function. 
	
	\begin{acknowledgments}
		The author gratefully acknowledges the Ministry of Education (Government of India) for the financial support via GATE fellowship (Registration No. MA19S43047204). 
	\end{acknowledgments}
	
	\section*{Data Availability}
	The data that support the plots within this paper and
	other findings of this study are  available
	within the article.
	\nocite{*}
	\section*{REFERENCES}
	\bibliography{aipsamp}
	
\end{document}